\documentclass{article}

\usepackage{amssymb}

\newcommand{\real}{{\mathbb R}}
\newcommand{\iin}{\!\in\!}

\newcommand{\Lip}{\mbox{\rm Lip}}
\newcommand{\Ub}{{\bf U_b}}
\newcommand{\Ubplus}{{\bf U^+_b}}
\newcommand{\Meas}{{\bf M}}
\newcommand{\UMeas}{{\bf M_u}}
\newcommand{\TMeasplus}{{\bf M^+_t}}
\newcommand{\UMeasplus}{{\bf M^+_u}}
\newcommand{\TMeas}{{\bf M_t}}
\newcommand{\Mol}{{\bf Mol}}
\newcommand{\Allmaps}{{\bf F}}
\newcommand{\sect}[1]{\setminus_{#1}}
\newcommand{\intx}[2]{#2_{\centerdot \, #1}}
\newcommand{\bigsep}{\mbox{\Large $|$}}
\newcommand{\conv}{\star}
\newcommand{\measprod}{\otimes}
\newcommand{\dirpr}{\!\!\times}
\newcommand{\semiu}{\raisebox{1mm}{$\ast$}}
\newcommand{\ru}{r\,}
\newcommand{\wstar}{weak\raisebox{1mm}{$\ast$}\ }
\newcommand{\compln}[1]{\widehat{#1}}
\newcommand{\compactn}[1]{\overline{#1}}
\newcommand{\qed}{\hfill $\Box$\vspace{3mm}}

\newtheorem{theorem}{Theorem}[section]
\newtheorem{lemma}[theorem]{Lemma}
\newtheorem{corollary}[theorem]{Corollary}

\title{Uniform measures and convolution on topological groups}

\author{Jan Pachl  \\ Toronto, Ontario, Canada}
\date{May 14, 2010 (version 4)}

\begin{document}
\maketitle

\begin{abstract}
Uniform measures are the functionals on the space of bounded uniformly continuous functions
that are continuous on every bounded uniformly equicontinuous set.
This paper describes the role of uniform measures in the study of
convolution on an arbitrary topological group.
\end{abstract}


\section{Introduction}
    \label{section:introduction}

Uniform continuity and equicontinuity properties have often been used in reasoning about
convolution of measures on groups and semigroups.
The aim of this paper is to point out the role of one such continuity property,
namely the one that characterizes so called {\em uniform measures\/},
in deriving the basic properties of convolution.

Uniform measures (under different names) were first studied by
Berezanski\v{\i}~\cite{Berezanskii1968},
Fedorova~\cite{Fedorova1967} and
LeCam~\cite{LeCam1970}~\cite{LeCam1970unp}.
The term {\em uniform measure\/} was used by Frol\'{\i}k~\cite{Frolik1973}.
Later developments are surveyed in publications by Deaibes~\cite{Deaibes1975}
and Frol\'{\i}k~\cite{Frolik1977}.

LeCam~\cite{LeCam1970unp} noted that the space of uniform measures
``arises naturally in various arguments about convolutions or Fourier transforms on linear spaces''.
LeCam's approach was later developed further by Caby~\cite{Caby1987},
who proved
several results about the continuity of convolution of uniform measures on abelian groups.

For general (not necessarily commutative) topological groups,
Csisz\'{a}r~\cite{Csiszar1971} defined a property equivalent to
being a positive uniform measure, which he called $\varrho$-continuity,
and used it to prove results about the continuity of convolution.
Additional results in this direction were proved by Tortrat~\cite{Tortrat1971}.
Recently Ferri and Neufang~\cite{Ferri-Neufang2006} described the topological centre
of the convolution algebra on certain topological groups, using uniform measures.

As noted by Csisz\'{a}r~\cite{Csiszar1970},
the natural way to define the convolution of two measures on a semigroup
is to take the image of the direct product measure under the semigroup operation.
That is the definition adopted in this paper, for measures on topological groups.
However,
as in the general approach described in~\cite{Hewitt-Ross1979} and~\cite{Pym1964},
the ``measures'' considered here are functionals on the space of bounded
uniformly continuous functions, not set functions on subsets of the underlying group.
Thus the direct product is defined as a functional on the space of
bounded uniformly continuous functions for a suitable uniformity on the product space.

The paper is organized as follows.
First we recall the definitions and basic results about uniform measures.
Next we derive basic properties of direct products of
functionals on the space of bounded uniformly continuous functions,
including Fubini's theorem for uniform measures.
Convolution on topological groups is then defined in terms of direct products.
Sections~\ref{section:centre} and~\ref{section:compactification}
include observations on topological centres and compactifications,
and the concluding section describes a generalization for structures more general than
topological groups.

The contributions of the present paper are
(1) the clarification of the role of uniform measures in the study of
convolution on topological groups;
and (2) the derivation of results about convolution from those about direct
product of measures,
where measures are rather general functionals.
As a result of (2), the approach immediately yields results for structures more general
than topological groups, as noted in the concluding section.

Version 1 of this paper was dated August 5, 2006.
This is version~4.
The following changes were made since version~1:
\begin{itemize}
\item New section~\ref{section:compactification}.
\item Additions and updates to the bibliography.
\item Minor edits and clarifications in the text.
\item Correction for Theorem~\ref{th:invmeans}.
\end{itemize}
Lemmas and theorems
in sections~\ref{section:uniformmeasures} to~\ref{section:centre}
have kept their numbers from version~1.


\section{Uniform measures}
    \label{section:uniformmeasures}

For simplicity, we consider all linear spaces to be over the field $\real$ of reals,
and a {\em function\/} is a mapping into~$\real$.
For any result in this paper it is an easy exercise to prove the corresponding
result with the field of complex numbers instead of $\real$.

All topological and uniform spaces are assumed to be Hausdorff.
We describe uniform spaces by uniformly continuous pseudometrics (\cite{Gillman-Jerison1960}, Chap.~15),
abbreviated as u.c.p.

If $Y$ is a set and $W$ is a uniform space, denote by $ \Allmaps ( Y, W ) $
the set of all mappings from $Y$ to $W$ and endow it with the uniformity
generated by the pseudometrics
\[
d \semiu ( f , f' )
= 1 \; \wedge \; \sup \; \{ \; d( f(y), f'(y) ) \; | \; y \in Y \; \} \;
, \;\;\; f, f' \in \Allmaps ( X, W ) \;
\]
where $d$ ranges over all u.c.p. on $W$.

Let $X$ and $W$ be two uniform spaces.
If $\cal P$ is a set of mappings from $X$ to $W$
and $d$ is a pseudometric on $W$,
define
\[
d_{\cal P} ( x, x' ) \; = 1 \; \wedge \; \sup \; \{ \; d( p(x), p(x') ) \; | \; p \in \cal P \; \}
\]
for $ x, x' \!\in\! X $.
Say that $\cal P$ is {\em uniformly equicontinuous\/}
if $ d_{\cal P} $ is a u.c.p on $X$ for each u.c.p. $d$ on $W$.
When $ {\cal P} = \{ p \} $ is a single-element set,
it is uniformly equicontinuous if and only if $p$ is uniformly continuous.

Note that $ \cal P $ is uniformly equicontinuous if and only if
the mapping $\chi_{\cal P} : X \rightarrow \Allmaps ( {\cal P} , W )$ defined by
$ \chi_{\cal P} (x) (p) =  p(x) $ for $ p \in \, {\cal P} $ is uniformly continuous.

When $d$ is a pseudometric on a set $X$, define
\[
\Lip ( d )  =  \{ f: X \rightarrow \real \; \bigsep \; |f(x)| \leq 1 \; \mbox{\rm and} \;
| f(x) - f(x') | \leq d(x,x') \; \mbox{\rm for all} \; x,x' \in X  \}
\]
Then $\Lip(d)$ is a compact subset of the product space $ \real^X $; we always consider $\Lip(d)$
with this compact topology.

When $X$ is a uniform space, denote by $\Ub(X)$ the space of bounded uniformly continuous
functions on $X$ with the norm
\[
\| f \| \; = \; \sup \{ \; | f(x) | \; \bigsep \; x \in X \} \;\; \mbox{\rm for} \;\; f \in \Ub(X)
\]
Denote by $\Meas(X)$ the norm dual of $\Ub(X)$; that is, the space of linear functions
$\mu$ on $\Ub(X)$ for which the norm
\[
\| \mu \| \; = \; \sup \{ \; | \mu(f) | \; \bigsep \; f \in \Ub(X) \;\; \mbox{\rm and} \;\;
\| f \| \leq 1 \; \}
\]
is finite.

When $X$ and $W$ are uniform spaces and $ p : X \rightarrow W $ is a
uniformly continuous mapping, define
$ \Meas(p) : \Meas(X) \rightarrow \Meas(W) $ by
\[
\Meas(p) ( \mu ) ( f ) = \mu ( f \circ p ) \;\; \mbox{\rm for}
\;\; \mu \in \Meas(X) \;\; \mbox{\rm and} \;\; f \in \Ub(W) \, .
\]
We often write simply $p(\mu)$ instead of $\Meas(p)(\mu)$.

Consider three subspaces of $ \Meas(X) $:
\begin{enumerate}
\item
For each $ x \in X $, define the {\em Dirac measure\/} $\delta_x \in \Meas(X)$ by
$ \delta_x (f) = f(x) $ for  $ f \in \Ub(X) $.
A linear function on $\Ub(X)$ is called a {\em molecular measure\/} if it is a linear combination
of finitely many Dirac measures.
In other words, it is of the form
\[
f \mapsto \sum_{x\in K} r(x) f(x) \;\;\; \mbox{\rm for} \;\;\; f \in \Ub(X)
\]
where $ K \subseteq X $ is finite and $ r(x) \in \real $ for $ x \in K $.
The space of molecular measures on $X$ is denoted $\Mol(X)$.
\item $ \TMeas(X) $ is the space of linear functions on $\Ub(X)$ that are continuous
on the unit $\| . \| $ ball in $\Ub(X)$
with respect to the compact-open topology.
That is, the linear functions on $\Ub(X)$ defined by integral with respect
to bounded tight (a.k.a. Radon) measures on $X$.
\item $ \UMeas(X) $ is the space of linear functions on $\Ub(X)$ that are continuous on
$\Lip(d)$ for every u.c.p. $d$ on $X$.
Here, as always, $\Lip(d)$ is considered with its compact topology,
as a subset of the product space $ \real^X $.
The elements of $ \UMeas(X) $ are called {\em uniform measures on~$X$\/}.
\end{enumerate}

Let $ \compactn{X} $ with a topological embedding $ \iota : X \rightarrow \compactn{X} $
be the uniform compactification of a uniform space $X$ (also called the Samuel compactification
--- \cite{Isbell1964}, II.32).
Then $ \Meas(\compactn{X}) = \TMeas(\compactn{X}) $ because $ \compactn{X} $ is compact,
and $ \Meas(\iota) $ is a Banach space isomorphism from
$ \Meas(X) $ to $ \Meas(\compactn{X}) $.
Thus the elements of $ \Meas(X) $ are in one to one correspondence with
bounded tight measures on $ \compactn{X} $.

Frol\'{\i}k~\cite{Frolik1973b} showed that the elements of $ \UMeas(X) $
may be usefully studied as measures on the compactification.
However, when considered on $X$ itself,
uniform measures need not be measures
in the sense of standard measure theory,
and in particular they need not be countably additive on $X$.
Nevertheless, the term {\em uniform measure\/} is established
in the literature, and we continue to use it in this paper.

\begin{theorem}
    \label{th:umeas1}
The following hold for any (Hausdorff) uniform space $X$.
\begin{enumerate}
\item
$ \Mol(X) \; \subseteq \; \TMeas(X) \; \subseteq \; \UMeas(X) \; \subseteq \; \Meas(X) \; . $
\item
If $ p : X \rightarrow W $ is a uniformly continuous mapping then $ \Meas(p) $
maps $ \Mol(X) $ into $ \Mol(W) $, $ \TMeas(X) $ into $ \TMeas(W) $, and
$ \UMeas(X) $ into $ \UMeas(W) $.
\item
    \label{th:umeas1:part3}
If $X$ is complete metric or uniformly locally compact then $ \TMeas(X) = \UMeas(X) $.
\item If $ \mu \in \Meas(X) $ then $ \mu \in \UMeas(X) $ if and only if
$ p (\mu) \in \TMeas(W) $ for every uniformly continuous mapping $ p : X \rightarrow W $
to a complete metric space $W$.
\end{enumerate}
\end{theorem}

\noindent
{\bf Proof.}
Parts~1 and~2 follow from the definition.
Part~3 for complete metric $X$ is proved in~\cite{Berezanskii1968},
\cite{Fedorova1967} and~\cite{LeCam1970unp},
and for uniformly locally compact $X$ in~\cite{Berezanskii1968} and~\cite{Fedorova1967}.
Part~4 is proved in~\cite{LeCam1970unp}.
\qed

In addition to the norm topology, we use two other topologies on $ \Meas(X) $:
\begin{itemize}
\item
The {\em \wstar topology\/} on $\Meas(X)$ is the weak topology
of the duality $ \langle \, \Meas(X), \Ub(X) \, \rangle $.
\item
The {\em UEB topology\/} on $\Meas(X)$ is the topology of uniform convergence
on the sets $\Lip(d)$, where $d$ ranges over all u.c.p. on $X$.
Equivalently, it is the topology of uniform convergence
on uniformly equicontinuous bounded subsets of $ \Ub(X) $.
The {\em UEB uniformity\/} is the corresponding translation-invariant
uniformity on $ \Meas(X) $.
\end{itemize}

Write
\begin{eqnarray*}
\Ubplus (X) \; & = & \; \{ \; f \in \Ub(X) \; | \; f(x) \geq 0 \; \mbox{\rm for} \;
x \in X \; \} \\
\Meas^+ (X) \; & = & \; \{ \; \mu \in \Meas(X) \; | \; \mu(f) \geq 0 \; \mbox{\rm for} \;
f \in \Ubplus (X) \; \}
\end{eqnarray*}
and define $\Mol^+ (X)$, $\TMeasplus (X)$ and $\UMeasplus (X)$ similarly.

Theorem~\ref{th:umeas-topology} enumerates basic results about topologies on $\UMeas(X)$.

\begin{theorem}
    \label{th:umeas-topology}
The following hold for any (Hausdorff) uniform space $X$.
\begin{enumerate}
\item \label{th:umeas-topology:part1}
The space $\Mol(X)$ is dense in $\UMeas(X)$ in the UEB topology.
\item \label{th:umeas-topology:part2}
The space $\Mol^+ (X)$ is dense in $\UMeasplus(X)$ in the UEB topology.
\item \label{th:umeas-topology:part3}
The space $\UMeas(X)$ with the UEB topology is complete.
\item \label{th:umeas-topology:part4}
The spaces $ \Meas(X) $ and $ \UMeas(X) $ with the $ \|.\| $ topology are complete.
\item \label{th:umeas-topology:part5}
The space $\UMeas(X)$ with the \wstar topology is sequentially complete.
\item \label{th:umeas-topology:part6}
Every \wstar compact subset of $\UMeas(X)$ is UEB compact.
\item \label{th:umeas-topology:part7}
If $ \mu \in \UMeasplus (X) $ then every UEB neighbourhood of $\mu$ in $\Meas^+(X)$
is also a \wstar neighbourhood of $\mu$.
\item \label{th:umeas-topology:part8}
The \wstar topology and the UEB topology coincide on $\UMeasplus (X)$.
\end{enumerate}
\end{theorem}

\noindent
{\bf Proof.}
As is pointed out in~\cite{Fedorova1967}, parts~\ref{th:umeas-topology:part1}
and~\ref{th:umeas-topology:part3} follow immediately from Grothendieck's completeness theorem
(\cite{Schaefer1966}, IV.6).
Part~\ref{th:umeas-topology:part2} is proved in~\cite{Berezanskii1968},
\cite{Fedorova1967} and~\cite{LeCam1970unp}.
In part~\ref{th:umeas-topology:part4},
the norm completeness of $ \Meas(X) $ follows from the definition,
and that of $ \UMeas(X) $ from part~\ref{th:umeas-topology:part3}.
Parts~\ref{th:umeas-topology:part5} and~\ref{th:umeas-topology:part6}
are proved in~\cite{Cooper-Schachermayer1981} and~\cite{Pachl1979}.
Parts~\ref{th:umeas-topology:part7} and~\ref{th:umeas-topology:part8}
are proved in~\cite{LeCam1970unp} (note that part~\ref{th:umeas-topology:part8} is
an immediate corollary of part~\ref{th:umeas-topology:part7}).
\qed


\section{Direct product}
    \label{section:directp}

To prepare the way for the definition of convolution in section~\ref{section:convolution},
in this section we define the direct product of functionals in $ \Meas(X) $ and $ \Meas(Y) $
on the semiuniform product of two uniform spaces $X$ and $Y$.

We adapt general lambda-calculus notation for the specific purposes of this paper,
to mark domain restriction for multivariate mappings.
If $p$ maps $x$ to $p(x)$ then $ \sect{x} p (x) $ is simply $p$ itself.
If $p$ maps $ (x,y) \in X \times Y $ to $ p(x,y) \in W $ then
the mappings $ \sect{x} p(x,y) : X \rightarrow W $ and
$ \sect{y} p(x,y) : Y \rightarrow W $
(the {\em sections of\/} $p$ along the two coordinates) are defined by
\[
\sect{x} p(x,y_0 )  ( x_0 ) = \sect{y} p(x_0 ,y)  ( y_0 ) = p(x_0 ,y_0 )
\]
for $ x_0 \in X $, $ y_0 \in Y $, and similarly for a mapping from the product of three or more sets.

Let $X$ and $Y$ be two uniform spaces.
The {\em semiuniform product\/} $ X \semiu Y $ of $X$ and $Y$ is a uniform space on the set
$X \dirpr Y$.
For the definition and basic properties of $ X \semiu Y $, see~Ch.~III in Isbell's book~\cite{Isbell1964}.
(Note that uniformly equicontinuous sets of mappings are called
{\em equiuniformly continuous\/} in~\cite{Isbell1964}.)
In what follows we only need the following characterization of $ X \semiu Y $ (loc. cit., Ch.~III,~Th.~22).
For any uniform space $W$,
a mapping $ p : X \semiu Y \rightarrow W$ is uniformly continuous if and only if
\begin{description}
\item[{\rm(SU1)}] the set $ \{ \sect{x} p(x,y)  \; | \; y \!\in\! Y  \} $ of mappings from $X$ to $W$
is uniformly equicontinuous; and
\item[{\rm(SU2)}] for each $ x \in X $, the mapping $ \sect{y} p(x,y)  $ from $Y$ to $W$
is uniformly continuous.
\end{description}

\begin{lemma}
    \label{lemma:equiu}
Let $X$, $Y$ and $W$ be uniform spaces.
A set $\cal P$ of mappings from $ X \semiu Y $ to $W$ is uniformly equicontinuous if and only if
\begin{description}
\item[{\rm(SU1*)}] the set $ \{ \sect{x} p(x,y)  \; | \; p \!\in\! {\cal P}, \; y \!\in\! Y \} $
of mappings from $X$ to $W$ is uniformly equicontinuous; and
\item[{\rm(SU2*)}] for each $ x \in X $,
the set $ \{ \sect{y} p(x,y)  \; | \; p \!\in\! \cal P \} $
of mappings from $Y$ to $W$ is uniformly equicontinuous.
\end{description}
\end{lemma}

\noindent
{\bf Proof.}
Apply (SU1) and (SU2) with $ \Allmaps ({\cal P}, W ) $ in place of $W$.
\qed

When $X$ and $Y$ are uniform spaces,
$ f \in \Ub(X \semiu Y) $, $ \nu \in \Meas(Y) $, define
$ \intx{\nu}{f} = \sect{x} \nu ( \sect{y} f(x,y)  ) $.
In view of (SU2),
the function  $ \intx{\nu}{f} $ is well defined on $X$.

\begin{lemma}
\label{lemma:dproduct0}
Let $X$ and $Y$ be uniform spaces,
let $ \cal F $ be a uniformly equicontinuous subset of $ \Ub(X \semiu Y) $,
and $ B \subseteq \Meas(Y) $ such that
$ s = \sup \{ \, \| \nu \| \, \bigsep \, \nu \in B \} $ is finite.
Then $ \| \intx{\nu}{f} \| \leq s \| f \| $
for each $ f \in {\cal F} $, $ \nu \in B $, and
the set
$ \{ \, \intx{\nu}{f} \; | \; f \in {\cal F} \, , \; \nu \in B \, \} $
is uniformly equicontinuous on $X$.
\end{lemma}

\noindent
{\bf Proof.}
Since $ \| \sect{y} f(x,y)  \| \leq \| f \| $ for each $x\in X$, we have
\[
| \intx{\nu}{f} (x) | \leq \| \nu \| \, . \, \| \sect{y} f(x,y)  \|
\leq  s \| f \|
\]
which shows that each $\intx{\nu}{f}$ is bounded
and $ \| \intx{\nu}{f} \| \leq s \| f \| $.

By (SU1*), there is a u.c.p. $d$ on $X$ such that
$
| f(x,y) - f(x',y) | \; \leq \; d(x,x')
$
for $ f \in {\cal F} $, $x,x' \in X$, and $y\in Y$.
Therefore for $ \nu \in B $ we have
\begin{eqnarray*}
\lefteqn{| \intx{\nu}{f} (x) - \intx{\nu}{f} (x') |
    = | \nu ( \sect{y} f(x,y)  ) - \nu ( \sect{y} f(x',y) ) | } \\
& = & | \nu ( \sect{y} f(x,y)  - \sect{y} f(x',y)  ) |
    \leq  \| \nu \| \, . \, d(x,x') \leq s \, . \, d(x,x')
\end{eqnarray*}
which shows that the set
$ \{ \, \intx{\nu}{f} \; | \; f \in {\cal F} \, , \; \nu \in B \, \} $
is uniformly equicontinuous.
\qed

\begin{lemma}
\label{lemma:dproduct1}
Let $X$ and $Y$ be uniform spaces, $ \nu \!\in\! \Meas(Y) $ and
$ f \!\in\! \Ub(X \semiu Y) $.
Then $ \intx{\nu}{f} \in \Ub(X) $
and $ \| \intx{\nu}{f} \| \leq \| \nu \| \, . \, \| f \| $.
\end{lemma}

\noindent
{\bf Proof} follows immediately from Lemma~\ref{lemma:dproduct0}, with
$ {\cal F} = \{ f \} $ and $ B = \{ \nu \} $.
\qed

Let $X$ and $Y$ be uniform spaces and $ \mu \in \Meas(X)$, $ \nu \in \Meas(Y) $.
For $ f \in \Ub( X \semiu Y) $,
define
\[
\mu \measprod \nu ( f ) = \mu ( \intx{\nu}{f} ) = \mu ( \,\sect{x} \nu ( \,\sect{y} f(x,y)  \,) \,)
\]
By Lemma~\ref{lemma:dproduct1}, $ \mu \measprod \nu ( f ) $ is well defined for each
$ f \in \Ub( X \semiu Y) $.
Call $ \mu \measprod \nu $ the {\em direct product\/} of $\mu$ and $\nu$.

Note that if $ g \in \Ub(X) $, $ h \in \Ub(Y) $
and $ f(x,y) = g(x)h(y) $ for $ x \in X$, $ y \in Y $,
then $ f \in \Ub( X \semiu Y )$ and
$ \mu \measprod \nu (f) = \mu(g) \, . \,\nu(h) $.

\begin{lemma}
\label{lemma:dproduct2}
Let $X$ and $Y$ be uniform spaces and $ \mu \!\in\! \Meas(X)$, $ \nu \!\in\! \Meas(Y) $.
Then $ \mu \measprod \nu \in \Meas(X \semiu Y) $ and
$ \| \mu \measprod \nu \| \leq \| \mu \| \, . \, \| \nu \| $.
\end{lemma}

\noindent
{\bf Proof.}
The linearity of $ \mu \measprod \nu $ follows from the linearity
of $\mu$ and $\nu$.
By Lemma~\ref{lemma:dproduct1},
\[
\| \mu \measprod \nu ( f ) \| = \| \mu ( \intx{\nu}{f} ) \|
\leq \| \mu \| \, . \, \| \intx{\nu}{f} \| \leq \| \mu \| \, . \, \| \nu \| \, . \, \| f \|
\]
and therefore $ \| \mu \measprod \nu \| \leq \| \mu \| \, . \, \| \nu \| $.
\qed

\begin{lemma}
    \label{lemma:measprod-closure}
Let $X$ and $Y$ be uniform spaces and $ \mu \in \Meas(X)$, $ \nu \in \Meas(Y) $.
\begin{enumerate}
\item
If $ \mu \in \Meas^+ (X) $ and $ \nu \in \Meas^+ (Y) $ then $ \mu \measprod \nu \in \Meas^+ ( X \semiu Y ) $.
\item
If $ \mu \in \Mol(X) $ and $ \nu \in \Mol(Y) $ then $ \mu \measprod \nu \in \Mol( X \semiu Y ) $.
\end{enumerate}
\end{lemma}

\noindent
{\bf Proof.}
Part 1 follows directly from the definition of $\measprod$.

2. If
\begin{eqnarray*}
\mu(f) & = & \sum_{x\in K}  r(x) \, f(x) \;\;\; \mbox{\rm for} \;\;\; f \in \Ub(X) \\
\nu(f) & = & \sum_{y\in K'} r'(y)\, f(y) \;\;\; \mbox{\rm for} \;\;\; f \in \Ub(Y)
\end{eqnarray*}
then
\[
\mu \measprod \nu (f) = \sum_{x\in K} \sum_{y\in K'} r(x) r'(y) f(x,y)
\;\;\; \mbox{\rm for} \;\;\; f \in \Ub(X \semiu Y)
\]
\qed

\begin{lemma}
    \label{lemma:dproduct-linear}
Let $X$ and $Y$ be uniform spaces,
$ \mu_1 , \mu_2 \in \Meas(X) $,
$ \nu_1 , \nu_2 \in \Meas(Y) $,
and $ r \in \real $.
Then
\begin{eqnarray*}
& & ( r \mu_1 ) \measprod \nu_1 = \mu_1 \measprod ( r \nu_1 ) = r ( \mu_1 \measprod \nu_1 )     \\
& & ( \mu_1 + \mu_2 ) \measprod \nu_1 = ( \mu_1 \measprod \nu_1 ) + ( \mu_2 \measprod \nu_1 )   \\
& & \mu_1 \measprod ( \nu_1 + \nu_2 ) = ( \mu_1 \measprod \nu_1 ) + ( \mu_1 \measprod \nu_2 )
\end{eqnarray*}
\end{lemma}

\noindent
{\bf Proof} follows directly from the definition of $\measprod$
and the linearity of $ \mu_1$, $\mu_2$, $\nu_1$ and $\nu_2$.
\qed

\begin{lemma}
    \label{lemma:dproduct3}
Let $X$, $Y$ and $Z$ be uniform spaces,
$ \mu \in \Meas(X) $, $ \nu \in \Meas(Y) $, and $ \xi \in \Meas(Z) $.
Denote by $ \varphi $ the canonical bijection
from the set $ ( X \times Y ) \times Z $ onto the set $ X \times ( Y \times Z ) $
defined by $ \varphi ( ( ( x,y ), z ) ) = ( x, ( y , z ) ) $.
Then
\begin{enumerate}
\item
$ \varphi: ( X \semiu Y ) \semiu Z \rightarrow X \semiu ( Y \semiu Z ) $
is a uniform isomorphism.
\item
$ \varphi ( ( \mu \measprod \nu ) \measprod \xi ) = \mu \measprod ( \nu \measprod \xi ) $.
\end{enumerate}
\end{lemma}

\noindent
{\bf Proof.}
1. We will show that any mapping $ p: X \semiu ( Y \semiu Z ) \rightarrow W $ to a uniform space $W$
is uniformly continuous if and only if
the mapping $ p \circ \varphi $ is uniformly continuous on $ ( X \semiu  Y ) \semiu Z $.

By the definition of the semiuniform product and by Lemma~\ref{lemma:equiu},
$ p \circ \varphi $ is uniformly continuous on $ ( X \semiu  Y ) \semiu Z $ if and only if \\
\indent (a) the set $ \{ \; \sect{x} p\circ\varphi((x,y),z)  \; | \; y \in Y, \; z \in Z \; \} $
is uniformly equicontinuous on $X$; \\
\indent (b) for each $ x \in X $,
the set $ \{ \; \sect{y} p\circ\varphi((x,y),z)  \; | \; z \in Z \; \} $
is uniformly equicontinuous on $Y$; \\
\indent (c) for each $ x \in X $ and $y \in Y $,
the mapping $ \sect{z} p\circ\varphi((x,y),z)  $
is uniformly continuous on $Z$.

\noindent
On the other hand, $ p $ is uniformly continuous
on $ X \semiu ( Y \semiu Z ) $ if and only if \\
\indent (a') the set $ \{ \; \sect{x} p(x,(y,z))  \; | \; y \in Y, \; z \in Z \; \} $
is uniformly equicontinuous on $X$; \\
\indent (b') for each $ x \in X $,
the set $ \{ \; \sect{y} p(x,(y,z))  \; | \; z \in Z \; \} $
is uniformly equicontinuous on $Y$; \\
\indent (c') for each $ x \in X $ and $y \in Y $,
the mapping $ \sect{z} p(x,(y,z))  $
is uniformly continuous on $Z$.

\noindent
Clearly (a), (b) and (c) are equivalent to (a'), (b') and (c').

2. For any $ f \in \Ub ( X \times ( Y \times Z ) ) $ we have
\begin{eqnarray*}
\lefteqn{\mu \measprod ( \nu \measprod \xi ) ( f ) } \\
& = & \mu ( \sect{x} ( \nu \measprod \xi ) ( \sect{(y,z)} f ( x, ( y, z ) ) ) )
    = \mu ( \sect{x}  \nu ( \sect{y} \xi ( \sect{z} f ( x, ( y, z ) ) ) ) )  \\
& = & \mu ( \sect{x}  \nu ( \sect{y} \xi ( \sect{z} f \circ \varphi ( ( x, y ) , z ) )  ) )
    = \mu \measprod \nu ( \sect{(x,y)} \xi ( \sect{z} f \circ \varphi ( ( x, y ) , z ) ) ) \\
& = & ( \mu \measprod \nu ) \measprod \xi ( f \circ \varphi )
\end{eqnarray*}
which shows that
$ \varphi ( ( \mu \measprod \nu ) \measprod \xi ) = \mu \measprod ( \nu \measprod \xi ) $.
\qed

Now we consider the continuity of the direct product operation $\measprod$.
Uniform measures have an important role in these results.
Theorem~\ref{th:measprod-cont} deals with the joint continuity of $\measprod$
in the UEB topology.
Corollary~\ref{cor:measprod-cont}
and Theorem~\ref{th:measprod-cont-positive}
deal with the joint continuity in the \wstar topology.

\begin{theorem}
    \label{th:measprod-cont}
Let $X$ and $Y$ be uniform spaces,
let $ B \subseteq \Meas(Y) $ be a set bounded in the $\| . \|$ norm on $ \Meas(Y)$,
$ \mu_0 \in \UMeas(X) $, and $ \nu_0 \in B $.
When $ \Meas(X) $, $ \Meas(Y) $ and $ \Meas ( X \semiu Y ) $ are endowed with their UEB topologies,
the mapping $ ( \mu, \nu ) \mapsto \mu \measprod \nu $
from $ \Meas(X) \times B $ to $ \Meas( X \semiu Y ) $
is jointly continuous at $ ( \mu_0 , \nu_0 ) $.
\end{theorem}

\noindent
{\bf Proof.}
Take any UEB neighbourhood of $ \mu_0 \measprod \nu_0 $
in $ \Meas (X \semiu Y) $ of the form
\[
N = \{ \; \xi \in \Meas(X \semiu Y) \; \bigsep
    \; | \xi ( f ) - \mu_0 \measprod \nu_0 ( f ) | < 4 \varepsilon
    \;\; \mbox{\rm for} \;\; f \in {\cal F} \; \}
\]
where ${\cal F} \subseteq \Ub ( X \semiu Y ) $ is
bounded and uniformly equicontinuous and $ \varepsilon > 0 $.

By Lemma~\ref{lemma:dproduct0}, the set
$ \{ \, \intx{\nu}{f} \; | \; f \in {\cal F} \, , \; \nu \in B \, \} $
is bounded and uniformly equicontinuous on $X$.
Therefore
\[
N_1 = \{ \; \mu \in \Meas(X) \; |
    \; | \mu ( \intx{\nu}{f} ) - \mu_0 ( \intx{\nu}{f} ) | < \varepsilon
    \;\; \mbox{\rm for} \;\; f \in {\cal F}, \; \nu \in B \; \}
\]
is a UEB neighbourhood of $ \mu_0 $ in $ \Meas(X) $.

By Theorem~\ref{th:umeas-topology} (part~\ref{th:umeas-topology:part1}),
there is $ \mu_1 \in \Mol(X) \cap N_1 $.
Now $\mu_1$ is of the form
\[
\mu_1 ( g ) = \sum_{x\in K} r(x) g(x) \;\;\; \mbox{\rm for} \;\;\; g \in \Ub(X)
\]
where $ K \subseteq X $ is finite and $ r(x) \in \real $ for $ x \in K $.
Set $ r = \sum_{x\in K} | r(x) | $.
By (SU2*), the set $ \{ \sect{y}{f(x,y)} \; | \; x \in K, f \in {\cal F} \} $ is
uniformly equicontinuous on $Y$ and clearly it is bounded in the $ \|. \|$ norm.
Define
\[
N_2 = \{ \; \nu \in B \; \bigsep
    \; | \nu ( \sect{y} f(x,y)  ) - \nu_0 ( \sect{y} f(x,y)  ) | < \varepsilon / r
    \;\; \mbox{\rm for} \;\; x \in K, f \in {\cal F} \; \}
\]
if $ r > 0 $ and $ N_2 = B $ if $ r = 0 $.
Then $ N_2 $ is a UEB neighbourhood of $ \nu_0 $ in $B$.
Note that for $ \nu \in N_2 $ and $ f \in {\cal F} $ we have
\[
| \mu_1 ( \intx{\nu}{f} ) - \mu_1 ( \intx{\nu_0}{f} ) |
= \bigsep \sum_{x\in K} r(x) \intx{\nu}{f} (x) - \sum_{x\in K} r(x) \intx{\nu_0}{f} (x) \bigsep
< \varepsilon
\]
If $ \mu \in N_1 $, $ \nu \in N_2 $ and $ f \in {\cal F} $ then
\begin{eqnarray*}
\lefteqn{ \mu \measprod \nu (f) - \mu_0 \measprod \nu_0 (f) } \\
    & = & \mu ( \intx{\nu}{f} ) - \mu_0 ( \intx{\nu_0}{f} )  \\
    & = & \mu ( \intx{\nu}{f} ) - \mu_0 ( \intx{\nu}{f} )
    + \mu_0 ( \intx{\nu}{f} ) - \mu_1 ( \intx{\nu}{f} )
    + \mu_1 ( \intx{\nu_0}{f} ) - \mu_0 ( \intx{\nu_0}{f} )
    + \mu_1 ( \intx{\nu}{f} ) - \mu_1 ( \intx{\nu_0}{f} )
\end{eqnarray*}
and therefore
\begin{eqnarray*}
\lefteqn{ \bigsep \mu \measprod \nu (f) - \mu_0 \measprod \nu_0 (f) \bigsep } \\
& \leq & \bigsep \mu ( \intx{\nu}{f} ) - \mu_0 ( \intx{\nu}{f} ) \bigsep
        + \bigsep \mu_0 ( \intx{\nu}{f} ) - \mu_1 ( \intx{\nu}{f} ) \bigsep
        + \bigsep \mu_1 ( \intx{\nu_0}{f} ) - \mu_0 ( \intx{\nu_0}{f} ) \bigsep
        + \bigsep \mu_1 ( \intx{\nu}{f} ) - \mu_1 ( \intx{\nu_0}{f} ) \bigsep \\
& < & 4 \, \varepsilon \, .
\end{eqnarray*}

Thus $ \mu \measprod \nu \in N $ whenever $ \mu \in N_1 $ and
$ \nu \in N_2 $.
\qed

\begin{corollary}
    \label{cor:measprod-cont}
Let $X$ and $Y$ be uniform spaces.
When $ \UMeas (X) $, $ \UMeas (Y) $ and $ \Meas ( X \semiu Y ) $
are endowed with their \wstar topologies,
the mapping $ ( \mu, \nu ) \mapsto \mu \measprod \nu $
is jointly sequentially continuous from
$ \UMeas (X) \times \UMeas (Y) $ to $ \Meas ( X \semiu Y ) $.
\end{corollary}

Note that, by Theorem~\ref{th:measprod-closure} below, the mapping
$ ( \mu , \nu ) \mapsto \mu \measprod \nu $ maps
$ \UMeas (X) \times \UMeas (Y) $ into $ \UMeas ( X \semiu Y ) $.

\noindent
{\bf Proof.}
Let $ \mu_n \in \UMeas (X) $, $ \nu_n \in \UMeas (Y) $ for $ n = 1,2, \ldots $,
$ \mu \in \UMeas (X) $, $ \nu \in \UMeas (Y) $
and
\[
\lim_n \mu_n = \mu \;\; , \;\;\;  \lim_n \nu_n = \nu
\]
in the \wstar topologies.
The set $ \{ \; \nu_n \, | \, n = 1,2, \ldots \; \} $
is $ \| . \| $ bounded
by the Banach-Steinhaus theorem (\cite{Schaefer1966}, III.4).
By Theorem~\ref{th:umeas-topology} (part~\ref{th:umeas-topology:part6}),
the UEB topology and the \wstar topology on $ \UMeas (X) $ and $ \UMeas (Y) $
have the same compact sets and therefore the same convergent sequences.
Thus
\[
\lim_n \mu_n = \mu \;\; , \;\;\;  \lim_n \nu_n = \nu
\]
in the UEB topologies, and Theorem~\ref{th:measprod-cont} shows that
\[
\lim_n \mu_n \measprod \nu_n = \mu \measprod \nu
\]
in the UEB topology and therefore also in the \wstar topology.
\qed

The next result is a direct product version of Csisz\'{a}r's Theorem~1
in~\cite{Csiszar1971}.

\begin{theorem}
    \label{th:measprod-cont-positive}
Let $X$ and $Y$ be uniform spaces,
$ \mu_0 \in \UMeasplus(X) $, and $ \nu_0 \in \Meas^+ (Y) $.
When $ \Meas^+(X) $, $ \Meas^+(Y) $ and $ \Meas^+ ( X \semiu Y ) $
are endowed with their \wstar topologies,
the mapping $ ( \mu, \nu ) \mapsto \mu \measprod \nu $
from $ \Meas^+(X) \times \Meas^+(Y) $ to $ \Meas^+ ( X \semiu Y ) $
is jointly continuous at $ ( \mu_0 , \nu_0 ) $.
\end{theorem}

\noindent
{\bf Proof} is almost identical to the proof of
Theorem~\ref{th:measprod-cont} above.
Take any \wstar neighbourhood of $ \mu_0 \measprod \nu_0 $
in $ \Meas^+ (X \semiu Y) $ of the form
\[
N = \{ \; \xi \in \Meas^+ (X \semiu Y) \; \bigsep
    \; | \xi ( f ) - \mu_0 \measprod \nu_0 ( f ) | < 4 \varepsilon
    \;\; \mbox{\rm for} \;\; f \in {\cal F} \; \}
\]
where ${\cal F} \subseteq \Ub ( X\semiu Y ) $ is finite and $ \varepsilon > 0 $.

Let $B = \{ \, \nu \in \Meas^+ (Y) \; | \; \nu (1) \leq \nu_0 (1) + 1 \, \} $.
Then $B$ is a \wstar neighbourhood of $\nu_0$ in $ \Meas^+ (Y) $.
By Lemma~\ref{lemma:dproduct0}, the set
$ \{ \, \intx{\nu}{f} \; | \; f \in {\cal F} \, , \; \nu \in B \, \} $
is bounded and uniformly equicontinuous on $X$.
Therefore
\[
N_1 = \{ \; \mu \in \Meas^+ (X) \; |
    \; | \mu ( \intx{\nu}{f} ) - \mu_0 ( \intx{\nu}{f} ) | < \varepsilon
    \;\; \mbox{\rm for} \;\; f \in {\cal F}, \; \nu \in B \; \}
\]
is a UEB neighbourhood of $ \mu_0 $ in $ \Meas^+ (X) $.

By Theorem~\ref{th:umeas-topology} (part~\ref{th:umeas-topology:part2}),
there is $ \mu_1 \in \Mol^+ (X) \cap N_1 $.
Now $\mu_1$ is of the form
\[
\mu_1 ( g ) = \sum_{x\in K} r(x) g(x) \;\;\; \mbox{\rm for} \;\;\; g \in \Ub(X)
\]
where $ K \subseteq X $ is finite and $ r(x) \in \real $, $ r(x) \geq0 $ for $ x \in K $.
Set $ r = \sum_{x\in K}  r(x) $, and define
\[
N_2 = \{ \; \nu \in B \; \bigsep
    \; | \nu ( \sect{y} f(x,y)  ) - \nu_0 ( \sect{y} f(x,y)  ) | < \varepsilon / r
    \;\; \mbox{\rm for} \;\; x \in K, f \in {\cal F} \; \}
\]
if $ r > 0 $ and $ N_2 = B $ if $ r = 0 $.
Then $ N_2 $ is a \wstar neighbourhood of $ \nu_0 $ in $ \Meas^+ (Y) $.
Note that for $ \nu \in N_2 $ and $ f \in {\cal F} $ we have
\[
| \mu_1 ( \intx{\nu}{f} ) - \mu_1 ( \intx{\nu_0}{f} ) |
= \bigsep \sum_{x\in K} r(x) \intx{\nu}{f} (x) - \sum_{x\in K} r(x) \intx{\nu_0}{f} (x) \bigsep
< \varepsilon
\]
Now by Theorem~\ref{th:umeas-topology} (part~\ref{th:umeas-topology:part7}),
$N_1$ is also a \wstar neighbourhood of $ \mu_0 $ in $ \Meas^+ (X) $.
If $ \mu \in N_1 $, $ \nu \in N_2 $ and $ f \in {\cal F} $ then
\begin{eqnarray*}
\lefteqn{ \mu \measprod \nu (f) - \mu_0 \measprod \nu_0 (f) } \\
    & = & \mu ( \intx{\nu}{f} ) - \mu_0 ( \intx{\nu_0}{f} )  \\
    & = & \mu ( \intx{\nu}{f} ) - \mu_0 ( \intx{\nu}{f} )
    + \mu_0 ( \intx{\nu}{f} ) - \mu_1 ( \intx{\nu}{f} )
    + \mu_1 ( \intx{\nu_0}{f} ) - \mu_0 ( \intx{\nu_0}{f} )
    + \mu_1 ( \intx{\nu}{f} ) - \mu_1 ( \intx{\nu_0}{f} )
\end{eqnarray*}
and therefore
\begin{eqnarray*}
\lefteqn{ \bigsep \mu \measprod \nu (f) - \mu_0 \measprod \nu_0 (f) \bigsep } \\
& \leq & \bigsep \mu ( \intx{\nu}{f} ) - \mu_0 ( \intx{\nu}{f} ) \bigsep
        + \bigsep \mu_0 ( \intx{\nu}{f} ) - \mu_1 ( \intx{\nu}{f} ) \bigsep
        + \bigsep \mu_1 ( \intx{\nu_0}{f} ) - \mu_0 ( \intx{\nu_0}{f} ) \bigsep
        + \bigsep \mu_1 ( \intx{\nu}{f} ) - \mu_1 ( \intx{\nu_0}{f} ) \bigsep \\
& < & 4 \, \varepsilon \, .
\end{eqnarray*}

Thus $ \mu \measprod \nu \in N $ whenever $ \mu \in N_1 $ and
$ \nu \in N_2 $.
\qed

We have defined $ \mu \measprod \nu $ by applying first $\nu$ and then $\mu$.
In the classical setting of countably additive measures,
Fubini's theorem (\cite{Hewitt-Ross1979}, 13.8) states that the order may be reversed:
\[
\int\int f(x,y) \ \mbox{\rm d}\mu(x) \ \mbox{\rm d}\nu(y) =
\int\int f(x,y) \ \mbox{\rm d}\nu(y) \ \mbox{\rm d}\mu(x)
\]
In our setting of measures as functionals on uniformly continuous functions,
Fubini's formula becomes
\[
\mu ( \; \sect{x} \nu (\, \sect{y} f(x,y) \, ) \; ) =
\nu ( \; \sect{y} \mu (\, \sect{x} f(x,y) \, ) \; )
\]
The next theorem states that this is true when $ \mu \in \UMeas(X)$.

\begin{theorem}
    \label{th:Fubini}
{\bf(Fubini's theorem)}
Let $X$ and $Y$ be uniform spaces, $ \mu \!\in\! \UMeas(X)$, $ \nu \!\in\! \Meas(Y) $,
and $ f \in \Ub(X \semiu Y) $.
Then $ \sect{y} \mu( \sect{x} f(x,y)) \in \Ub(Y) $ and
$ \mu \measprod \nu \,( f ) = \nu ( \sect{y} \mu( \sect{x} f(x,y) ) ) $.
\end{theorem}

\noindent
{\bf Proof.}
First observe that the theorem holds when $ \mu \in \Mol(X) $, by the linearity of $\nu$.
In the general case we approximate $ \mu $ by molecular measures, as follows.

By Theorem~\ref{th:umeas-topology} (part~\ref{th:umeas-topology:part1}),
there is a net $ \{ \mu_\alpha \}_\alpha $ such that $ \mu_\alpha \in \Mol(X) $ for each $\alpha$
and $ \lim_\alpha \mu_\alpha = \mu $ in the UEB topology on $ \UMeas(X) $.

Take any $ f \in \Ub(X \semiu Y) $.
The set $ \{ \sect{x}{f(x,y)} \; | \; y \!\in\! Y  \} $
is uniformly equicontinuous on $X$ by (SU1), and  $\| . \| $ bounded.
Therefore
$
\lim_\alpha \; \sect{y} \mu_\alpha ( \sect{x}{f(x,y)} ) = \sect{y} \mu ( \sect{x}{f(x,y)} )
$
uniformly on $ Y $.

It follows that $ \sect{y} \mu ( \sect{x}{f(x,y)} ) \in \Ub(Y) $
and $ \lim_\alpha \nu ( \sect{y} \mu_\alpha ( \sect{x}{f(x,y)} ) )
      = \nu ( \sect{y} \mu ( \sect{x}{f(x,y)} ) ) $,
and from Theorem~\ref{th:measprod-cont} we get
\[
\mu \measprod \nu (f) = \lim_\alpha \mu_\alpha \measprod \nu (f)
= \lim_\alpha \nu ( \sect{y} \mu_\alpha ( \sect{x}{f(x,y)} ) )
= \nu ( \sect{y} \mu ( \sect{x}{f(x,y)} ) )
\]
where the second equality follows from the observation at the beginning of the proof.
\qed

Next we show that the property in Fubini's theorem as stated above holds only for uniform measures,
and is equivalent to a certain continuity property of the $\measprod$ operation.

\begin{theorem}
    \label{th:Fubiniconverse}
Let $X$ be a uniform space.
The following properties of $ \mu \in \Meas (X) $ are equivalent:
\begin{description}
\item[{\it (i)}]
$ \mu \in \UMeas (X) $.
\item[{\it (ii)}]
For every uniform space $Y$, if $ \;\nu \in \Meas(Y) $ and $ f \in \Ub ( X \semiu Y ) $
then $ \sect{y} \mu( \sect{x} f(x,y)) \in \Ub(Y) $ and
\[
\mu \measprod \nu \,( f ) = \nu ( \sect{y} \mu( \sect{x} f(x,y) ) ) \; .
\]
\item[{\it (iii)}]
For every uniform space $Y$, the mapping
$ \nu \mapsto \mu \measprod \nu $ from $ \Meas(Y) $
to $ \Meas ( X \semiu Y ) $ is continuous with respect to the
\wstar topologies on $ \Meas(Y) $ and $ \Meas ( X \semiu Y ) $.
\item[{\it (iv)}]
For every uniform space $Y$, the mapping
$ y \mapsto \mu \measprod \delta_y $ from $Y$
to $ \Meas ( X \semiu Y ) $ is continuous with respect to the
\wstar topology on $ \Meas ( X \semiu Y ) $.
\end{description}
\end{theorem}

\noindent
{\bf Proof.} (i) implies (ii) by Theorem~\ref{th:Fubini}.

Now assume (ii), and take any
$ \nu \in \Meas (Y) $ and a net $ \{ \nu_\alpha \}_\alpha $
such that $ \nu_\alpha \in \Meas(Y) $ for each $\alpha$
and $ \lim_\alpha \nu_\alpha = \nu $ in the \wstar topology on $ \Meas (Y)$.
Thus $ \lim_\alpha \nu_\alpha (h) = \nu (h) $ for each  $ h \in \Ub (Y) $.
If $ f \in \Ub ( X \semiu Y ) $ then by (ii) we get
$ \sect{y} \mu( \sect{x} f(x,y)) \in \Ub(Y) $
and
\[
\lim_\alpha \mu \measprod \nu_\alpha (f)
= \lim_\alpha \nu_\alpha (\sect{y} \mu( \sect{x} f(x,y)))
= \nu (\sect{y} \mu( \sect{x} f(x,y)) )
= \mu \measprod \nu (f) \;.
\]
That proves that (ii) implies (iii).

To prove that (iii) implies (iv), note that if $ y_\alpha \rightarrow y $
in the topology of $Y$ then $ \delta_{y_\alpha} \rightarrow \delta_y$
in the \wstar topology on $\Meas(Y)$.

To prove that (iv) implies (i),
let $d$ be a u.c.p. on $X$ and let $ \{ g_\alpha \}_{\alpha \in A} $
be a net such that $ g_\alpha \in \Lip(d) $ for each $ \alpha \in A $
and $ \lim_\alpha g_\alpha (x) = 0 $ for each $ x \in X $.
We wish to prove that $ \lim_\alpha \mu ( g_\alpha ) = 0 $, assuming (iv).
Choose an element $ \infty $ not in $A$ and
define a uniform space $Y$ on the set $ A \cup \{ \infty \} $ so that
\[
\Ub (Y) = \{ \; f : Y \rightarrow \real \;\; | \;\; f \; \mbox{\rm is bounded and}
    \; \lim_\alpha f ( \alpha ) = f (\infty) \; \} \; .
\]
Then the function defined by $ f(x, \alpha) = g_\alpha (x) $
and $ f(x,\infty) = 0 $ belongs to $ \Ub ( X \semiu Y ) $.
We have $ \lim \alpha = \infty $ in $Y$
and $ \sect{x} \delta_\alpha ( \sect{y} f(x,y) ) = g_\alpha $,
hence
\[
\lim_\alpha \mu ( g_\alpha ) =
\lim_\alpha \mu \measprod \delta_\alpha (f) = \mu \measprod \delta_\infty (f) = 0
\]
where the second equality follows from (iv).
\qed

To conclude this section, we show that
$ {\bf M}_u $ is closed under the direct product operation $\measprod$.

\begin{theorem}
    \label{th:measprod-closure}
Let $X$ and $Y$ be uniform spaces.
If $ \mu \in \UMeas(X) $ and $ \nu \in \UMeas(Y) $ then $ \mu \measprod \nu \in \UMeas( X \semiu Y ) $.
\end{theorem}

\noindent
{\bf Proof.}
First consider the case of $ \nu \in \UMeas (Y) $
and $ \mu = \delta_x $ (Dirac measure) for some $ x \in X $.
In this case we have
\[
\mu \measprod \nu ( f ) = \nu ( \sect{y} f(x,y) ) \;\;\; \mbox{\rm for} \;\;\; f \in \Ub(X \semiu Y)
\]
and therefore $ \mu \measprod \nu \in \UMeas ( X \semiu Y ) $.
By linearity, the same is true when $ \mu \in \Mol(X) $.

Now in the general case of $ \mu \in \UMeas(X) $ and $ \nu \in \UMeas(Y) $,
we approximate $ \mu $ by molecular measures,
as in the proof of Theorem~\ref{th:Fubini}.
By Theorem~\ref{th:umeas-topology} (part~\ref{th:umeas-topology:part1}),
there is a net $ \{ \mu_\alpha \}_\alpha $ such that $ \mu_\alpha \in \Mol(X) $ for each $\alpha$
and $ \lim_\alpha \mu_\alpha = \mu $ in the UEB topology on $ \UMeas(X) $.

By the argument above, $ \mu_\alpha \measprod \nu \in \UMeas ( X \semiu Y ) $ for each $\alpha$,
and by Theorem~\ref{th:measprod-cont} we have
$
\lim_\alpha \mu_\alpha \measprod \nu = \mu \measprod \nu
$
in the UEB topology on $ \Meas ( X \semiu Y ) $.
Therefore $ \mu \measprod \nu \in \UMeas(X \semiu Y)$
by Theorem~\ref{th:umeas-topology} (part~\ref{th:umeas-topology:part3}).
\qed


\section{Convolution on topological groups}
    \label{section:convolution}

Let $G$ be a group and $d$ a pseudometric on $G$.
Say that $d$ is {\em right-invariant\/}
if $ d(x,y) = d(xz,yz) $ for all $x,y,z \in G $.

When $G$ is a topological group,
the {\em right uniformity\/} is generated by all right-invariant
continuous pseudometrics on $G$.
We denote $ \ru G $ the set $G$ with the right uniformity.
If $G$ and $G'$ are topological groups and
$ h : G \rightarrow G' $ is a continuous homomorphism then
$h$ is uniformly continuous from $ \ru G $ to $ \ru G' $.

{\bf Warning.} The terminology and notation in this paper are based on those
commonly used in topology and measure theory.
Publications on abstract harmonic analysis often use another terminology and
notation, in which $ M(G) $ is the space of bounded tight measures on $G$,
$ \Ub (\ru G) $ is denoted LUC$(G) $, and the elements of LUC$(G) $ are called
{\em left\/} uniformly continuous.
\begin{center}
\begin{tabular}{|l|l|} \hline
{in this paper} \hspace{1.5cm} & {alternative notation} \\ \hline
$ \Ub(\ru G) $ & LUC$(G) $ \\ \hline
$ \Meas(\ru G) $ & LUC$(G)\raisebox{1mm}{$\ast$} $ \\ \hline
$ \TMeas(\ru G) $  & $ M(G) $ \\ \hline
$ \UMeas(\ru G) $  &  Leb$(G)\;\;$ \cite{Ferri-Neufang2006} \\ \hline
\end{tabular}
\end{center}

When $G$ is a topological group, denote by $m_G$ the binary group operation in $G$;
that is, the mapping $ m_G : G \times G \rightarrow G $ defined by $ m_G (x,y) = x\, y $.
The following lemma is a translation of well known properties of topological groups
to the language of semidirect products.

\begin{lemma}
    \label{lemma:rightu}
If $G$ is a topological group, then the group operation $ m_G $
is uniformly continuous from $ \ru G \semiu \ru G $ to $ \ru G $.
\end{lemma}

\noindent
{\bf Proof.}
We need to verify (SU1) and (SU2) in section~\ref{section:directp} with $ p = m_G $,
$ X = Y = W = \ru G $.

(SU1) follows directly from the definition, since for every right-invariant pseudometric $d$ we have
\[
d ( m_G ( x , y ) , m_G ( x' , y ) ) = d ( xy, x'y ) = d ( x , x' ) \; .
\]
To verify (SU2), fix $x \in G$ and a right-invariant pseudometric $d$ continuous on $G$.
Then the pseudometric $d_{\sect{y} m_G(x,y)}$ defined by
\[
d_{\sect{y}{m_G(x,y)}} (y,y') = d ( m_G ( x , y ) , m_G ( x , y' ) ) = d ( xy, xy' )
\]
is right-invariant and continuous, and therefore the mapping $ \sect{y}{m_G(x,y)} $
is uniformly continuous from $\ru G$ to $\ru G$.
\qed

Now let $G$ be a topological group, and
$ \mu, \nu \in \Meas(\ru G) $.
In section~\ref{section:directp} we defined the direct product
$ \mu \measprod \nu \in \Meas( \ru G \semiu \ru G ) $.
By Lemma~\ref{lemma:rightu}, $ \Meas(m_G) ( \mu \measprod \nu ) $ is well defined as
an element of $ \Meas(G)$.
Denote
\[
\mu \conv \nu = \Meas(m_G) ( \mu \measprod \nu )
\]
and call $ \mu \conv \nu $ the {\em convolution\/} of $\mu$ and $\nu$.

Thus we have
$
\mu \conv \nu ( f ) = \mu ( \sect{x} \nu ( \sect{y} f ( xy ) ) )
$
for $ \mu, \nu \in \Meas(\ru G) $ and $ f \in \Ub (\ru G) $.
It is easy to check that this definition agrees with the general definition
of convolution in~\cite{Hewitt-Ross1979} (19.1; see also 19.23(b)) and
with the general definition of {\em evolution\/}
in~\cite{Pym1964}.

If $ \mu \in \UMeas (\ru G) $ then also
\[
\mu \conv \nu ( f ) = \nu ( \sect{y} \mu ( \sect{x} f ( xy ) ) )
\]
by Fubini's Theorem~\ref{th:Fubini}.
In the terminology of Pym~\cite{Pym1964},
this means that the convolution and the evolution of $\mu$ and $\nu$
coincide when $ \mu \in \UMeas (\ru G) $.

The results in section~\ref{section:directp}
now immediately yield the properties of convolution in the following three theorems,
showing that $ \Meas(\ru G) $ and $ \UMeas(\ru G) $ with
the operations $ \conv $ and $ + $ and the norm $ \| . \| $ are Banach algebras.
For $ \Meas(\ru G) $ these results are well known~\cite{Hewitt-Ross1979}, \cite{Pym1965}.

\begin{theorem}
    \label{th:convolution-identities}
Let $G$ be a topological group,
$ \mu , \nu, \xi \in \Meas(\ru G) $,
and $ r \in \real $.
Then
\begin{eqnarray*}
& & ( \mu \conv \nu ) \conv \xi = \mu \conv ( \nu \conv \xi ) \\
& & ( r \mu ) \conv \nu = \mu \conv ( r \nu ) = r ( \mu \conv \nu )     \\
& & ( \mu + \nu ) \conv \xi = ( \mu \conv \xi ) + ( \nu \conv \xi )   \\
& & \mu \conv ( \nu + \xi ) = ( \mu \conv \nu ) + ( \mu \conv \xi )
\end{eqnarray*}
\end{theorem}

\noindent
{\bf Proof.} The first identity follows from Lemma~\ref{lemma:dproduct3}.
The remaining three identities follow from Lemma~\ref{lemma:dproduct-linear}.
\qed

\begin{theorem}
    \label{th:conv-closure}
For every topological group $G$, the sets $\Meas^+ (\ru G)$,
$ \Mol(\ru G) $, $ \Mol^+ (\ru G) $,
$ \UMeas(\ru G) $ and $ \UMeasplus(\ru G) $
are closed under the operation $\conv$.
\end{theorem}

\noindent
{\bf Proof} follows from Lemma~\ref{lemma:measprod-closure} and
Theorem~\ref{th:measprod-closure}.
\qed

In the next theorem, part~\ref{th:convolution-topology:part1} is well known.
A weaker version of
part~\ref{th:convolution-topology:part2} for commutative groups
is Theorem~3.2 in~\cite{Caby1987},
part~\ref{th:convolution-topology:part4} is a variant of Theorem~1 in~\cite{Csiszar1971},
and part~\ref{th:convolution-topology:part5} is included in Proposition~4.2
in~\cite{Ferri-Neufang2006}.
Part~\ref{th:convolution-topology:part2} is also similar to Lemma~2.2
in~\cite{Pym1965}.

\begin{theorem}
    \label{th:convolution-topology}
The following hold in any topological group $G$.
\begin{enumerate}
\item
    \label{th:convolution-topology:part1}
$ \| \mu \conv \nu \| \leq \| \mu \| . \| \nu \| $ for any $ \mu , \nu \in \Meas(\ru G) $,
\item
    \label{th:convolution-topology:part2}
Let $ B \subseteq \Meas(\ru G) $ be a set bounded in the $\| . \|$ norm,
$ \mu_0 \in \UMeas(\ru G) $, and $ \nu_0 \in B $.
When $ \Meas(\ru G)$ and its subset $B$ are endowed with their UEB topology,
the mapping $ ( \mu, \nu ) \mapsto \mu \conv \nu $
from $ \Meas(\ru G) \times B $ to $ \Meas(\ru G) $
is jointly continuous at $ ( \mu_0 , \nu_0 ) $.
\item
    \label{th:convolution-topology:part3}
When $ \UMeas(\ru G) $ is endowed with its \wstar topology,
the mapping $ ( \mu, \nu ) \mapsto \mu \conv \nu $
from $ \UMeas(\ru G) \times \UMeas(\ru G) $ to $ \UMeas(\ru G) $
is jointly sequentially continuous.
\item
    \label{th:convolution-topology:part4}
Let $ \mu_0 \in \UMeasplus(\ru G) $ and $ \nu_0 \in \Meas^+ (\ru G) $.
When $ \Meas^+(\ru G) $ is endowed with its \wstar topology,
the mapping $ ( \mu, \nu ) \mapsto \mu \conv \nu $
from $ \Meas^+(\ru G) \times \Meas^+(\ru G) $ to $ \Meas^+ ( \ru G ) $
is jointly continuous at $ ( \mu_0 , \nu_0 ) $.
\item
    \label{th:convolution-topology:part5}
If $ \mu \in \UMeas(\ru G) $ then the mapping $ \nu \mapsto \mu \conv \nu $
from $ \Meas(\ru G) $ to itself is \wstar continuous.
\end{enumerate}
\end{theorem}

\noindent
{\bf Proof.}
1. For any $ \mu , \nu \in \Meas(\ru G) $ and $ f \in \Ub (\ru G) $
we have $ \| f \| = \| f \circ m_G \| $,
and therefore
\[
| \mu \conv \nu ( f ) | \; = \; | \mu \measprod \nu ( f \circ m_G ) | \;
\leq \; \| \mu \measprod \nu \| \, . \, \| f \circ m_G \| \;
\leq \; \| \mu \| \, . \, \| \nu \| \, . \, \| f \| \;
\]
which shows that $ \| \mu \conv \nu \| \leq \| \mu \| . \| \nu \| $.

Part~\ref{th:convolution-topology:part2}
follows from Theorem~\ref{th:measprod-cont},
part~\ref{th:convolution-topology:part3}
from Corollary~\ref{cor:measprod-cont},
part~\ref{th:convolution-topology:part4}
from Theorem~\ref{th:measprod-cont-positive}
and part~\ref{th:convolution-topology:part5}
from Theorem~\ref{th:Fubiniconverse}.
\qed

For molecular measures, the definition of convolution yields an explicit formula:
If
\begin{eqnarray*}
\mu(f) & = & \sum_{x\in K}  r(x) \, f(x) \\
\nu(f) & = & \sum_{y\in K'} r'(y)\, f(y)
\end{eqnarray*}
for $ f \in \Ub(\ru G) $ and finite sets $K$, $K'$, then
\[
\mu \conv \nu (f) = \sum_{x\in K} \sum_{y\in K'} r(x) r'(y) f(xy) \; .
\]

Since $\Mol(\ru G)$ is dense in $\UMeas(\ru G)$ in the UEB topology,
part~\ref{th:convolution-topology:part2} in
Theorem~\ref{th:convolution-topology}
implies that algebraic identities
satisfied by the convolution operation on
$\Mol(\ru G)$ are inherited by $\UMeas(\ru G)$.
In particular, if $G$ is commutative then so is $\UMeas(\ru G)$.
That was proved by LeCam for the additive groups of locally convex spaces
(\cite{LeCam1970}, Prop.~5); however, his proof is valid for all commutative groups.
On the other hand, as is pointed out in~\cite{Csiszar1971},
(\cite{Hewitt-Ross1979},~19.24) and (\cite{Pym1964},~3.8),
$\Meas(\ru G)$ need not be commutative when $G$ is.

When $G$ is locally compact, $ \UMeas( \ru G ) = \TMeas( \ru G ) $
by Theorem~\ref{th:umeas1} (part~\ref{th:umeas1:part3}).
Some properties of convolution in $\TMeas(\ru G)$ on locally compact groups $G$
generalize to $\UMeas(\ru G)$ on arbitrary topological groups, but some do not.
Csisz\'{a}r~\cite{Csiszar1971} proved the following property for
a class of topological groups that includes abelian and locally compact groups:
If $ \mu \in \UMeasplus (\ru G) $, $ \nu \in \Meas^+ (\ru G) $
and $ \mu \conv \nu \in \UMeasplus (\ru G) $ then  $ \nu \in \UMeasplus (\ru G) $
(loc. cit., Lemma~2 and Lemma~3).
The next theorem shows that there are topological groups that do not have this property.

\begin{theorem}
    \label{th:metrizableexample}
For any metrizable topological group $G$, these two conditions are
equivalent:
\begin{description}
\item[{\it (i)}]
$G$ has a group completion.
\item[{\it (ii)}]
If $ \mu \in \UMeasplus (\ru G) $, $ \nu \in \Meas^+ (\ru G) $
and $ \mu \conv \nu \in \UMeasplus (\ru G) $ then  $ \nu \in \UMeasplus (\ru G) $.
\end{description}
\end{theorem}

\noindent
{\bf Proof.}
To prove that (i) implies (ii), assume that $G$ has a group completion.
Thus $G$ is a dense subgroup of a topological group $G'$
such that the uniform space $ \ru G' $ is a completion
of the uniform space $ \ru G $.
There is a natural isomorphism between $ \Ub(\ru G) $ and $ \Ub(\ru G') $,
between $ \Meas(\ru G) $ and $ \Meas(\ru G') $,
and between $ \UMeas(\ru G) $ and $ \UMeas(\ru G') $.
Since $\ru G'$ is a complete metric space, $ \UMeas(\ru G') = \TMeas(\ru G')$
by Theorem~\ref{th:umeas1} (part~\ref{th:umeas1:part3}),
and the statement (ii) follows from Lemma~4 in~\cite{Csiszar1971}.

To prove that (ii) implies (i), assume that $G$ does not have a group completion.
By Theorem~10.5 in~\cite{Roelcke-Dierolf1981},
there exists a set $ A \subseteq G $ such that $A$ is precompact in $ \ru G $
and $ \{ x^{-1} | x \in A \} $ is not.
Therefore there exist elements $ x_n \in G $, $ n =1,2,\ldots$,
such that the sequence $ \{ x_n \}_n $ is Cauchy in $ \ru G $
and the sequence $ \{ x^{-1}_n \}_n  $ is uniformly discrete in $ \ru G $.

Define $ \mu (f) = \lim_n f(x_n) $ for each $ f \in \Ub (\ru G)$.
Then $ \mu \in \UMeasplus(\ru G) $ by parts~\ref{th:umeas-topology:part3}
and~\ref{th:umeas-topology:part8}
in Theorem~\ref{th:umeas-topology} (or simply by observing that $\mu$
is the Dirac measure at the limit of $\{ x_n \}_n $ in the completion of $ \ru G $).
Take a free ultrafilter ${\cal U}$ on $ \{1,2,\ldots\} $ and define
$ \nu \in \Meas^+(\ru G) $ by
$ \nu (f) = \lim_{\cal U} f(x^{-1}_n) $ for each $ f \in \Ub (\ru G)$.

By Theorem~\ref{th:convolution-topology}
(part~\ref{th:convolution-topology:part2}
or~\ref{th:convolution-topology:part4}),
$
\mu \conv \nu = \lim_{\cal U} \;\;\delta_{x_n} \conv \;\delta_{x^{-1}_n} =
\delta_e
$
where $e$ is the unit of $G$, and therefore
$ \mu \conv \nu \in \Mol^+ (\ru G) \subseteq \UMeasplus (\ru G) $.

Since the sequence $ \{ x^{-1}_n \}_n  $ is uniformly discrete in $ \ru G $,
there is a right-invariant continuous pseudometric $d$ on $G$ such that
$ d( x^{-1}_m , x^{-1}_n ) > 2 $ for $ m \neq n $.
Define the functions $ f_n , f $ by
$ f_n (x) = 1 \wedge \min \{ \; d( x , x^{-1}_i ) \; | \; i = 1,2,\ldots,n \; \} $
and $ f(x) = \inf \{ \; f_n (x) \; | \; n =1,2,\ldots \; \} $ for $ x \in G $.
The set $ \{ \; f_n \; | \; n =1,2,\ldots \; \} \subseteq \Ub (\ru G) $
is bounded and uniformly equicontinuous and
$ \lim_n f_n (x) = f(x) $ for each $ x \in G $.
However, $ \nu ( f_n ) = 1 $ for each $n$ and $ \nu ( f ) = 0 $,
which shows that $ \nu \not\in \UMeasplus (\ru G) $.

This concludes the proof that if $G$ does not have a group completion
then (ii) does not hold in $G$.
\qed

Examples of metrizable groups that do not have a group completion
can be found in~\cite{Roelcke-Dierolf1981}.
Thus Theorem~\ref{th:metrizableexample}
answers Csisz\'{a}r's question on page 36 of~\cite{Csiszar1971}.
In the notation of that paper,
if $X$ is a metrizable group that does not have a group completion then
$ {\cal M}^{\varrho}_r (X)$ cannot replace $ {\cal M}^{p}_r (X)$
in Lemma~2 of~\cite{Csiszar1971},
and the inclusion
$ {\cal M}^{\varrho}_r (X) \subseteq {\cal M}^{p}_r (X) $ does not hold.


\section{Topological centre and amenability}
    \label{section:centre}

Let $G$ be a topological group.
Define
\[
Z =  \{ \; \mu \in \Meas(\ru G) \;
| \;
\mbox{\rm the mapping} \;\; \nu \mapsto \mu \conv \nu  \;\;
\mbox{\rm is \wstar continuous on} \; \Meas(\ru G) \; \}
\]
and call $Z$ the {\em topological centre\/} of $ \Meas(\ru G) $.
In the notation of~\cite{BJM1989}, $ Z = \Lambda ( \Meas(\ru G ) )$.

For $ \mu \in \Meas(\ru G) $ we have $ \mu \in Z $ if and only if
$\nu ( \sect{y} \mu ( \sect{x} f ( xy ) ) ) $ is defined and
\[
\mu ( \sect{x} \nu ( \sect{y} f ( xy ) ) )
= \nu ( \sect{y} \mu ( \sect{x} f ( xy ) ) )
\]
for all $ \nu \in \Meas(\ru G)$, $f \in \Ub(\ru G)$
(\cite{Lau1986}, Lemma~2).

From part~\ref{th:convolution-topology:part5} of Theorem~\ref{th:convolution-topology}
(alternatively, from Fubini's Theorem~\ref{th:Fubini})
it follows that
$ \UMeas(\ru G) \subseteq Z $
--- proved as part of Prop~4.2 in~\cite{Ferri-Neufang2006}.
For topological groups, this is a generalization of Lemma~3.1(a) of Wong~\cite{Wong1972}.

In view of Theorem~\ref{th:Fubiniconverse},
it is reasonable to ask whether $ \UMeas(\ru G) = Z $ for every topological group $G$.
Lau~\cite{Lau1986} proved $ \TMeas(\ru G) = Z $ for all locally compact groups and,
as noted above, $ \UMeas(\ru G) = \TMeas(\ru G) $ when $G$ is locally compact.
Ferri and Neufang~\cite{Ferri-Neufang2006}
proved $ \UMeas(\ru G) = Z $
for all $\omega$-bounded (not necessarily locally compact) topological groups.

The groups for which $ \UMeas(\ru G) = Z $ have a number of interesting properties;
see the corollaries in~\cite{Lau1986}, section 4.
In Theorem~\ref{th:uniqmprecompact} we point out one of these properties, namely the connection to unique amenability.

For any $ \mu \in \Meas(\ru G)$ and $ x \in G $,
the left translation of $\mu$ by $x$
is $ \delta_x \conv \mu $ where $ \delta_x $ is the Dirac measure.
If $ \mu \in \Meas(\ru G) $ and $ \delta_x \conv \mu = \mu $ for each $ x \in G $
then $ \mu $ is called {\em left-invariant\/}.
A functional $ \mu \in \Meas^+ (\ru G) $ such that $ \mu (1) = 1 $ is called a {\em mean\/}.
If there exists a left-invariant mean in $ \Meas^+ (\ru G) $ then $G$ is called
{\em (left) amenable\/}.

\begin{lemma}
    \label{lemma:uniqmean}
If $G$ is a topological group and $ \mu \in \Meas^+ (\ru G) $ is a unique
left-invariant mean then $ \mu \in Z $.
\end{lemma}

\noindent
{\bf Proof:} See the proof of Corollary 5 in~\cite{Lau1986}.
\qed

A topological group that admits a unique left-invariant mean
is called {\em uniquely (left) amenable}.
Every precompact topological group is uniquely amenable.
Megrelishvili, Pestov and Uspenskij (\cite{MPU2001}, 3.5) ask whether conversely every
uniquely amenable topological group is precompact.

\begin{theorem}
    \label{th:uniqmprecompact}
If $G$ is a uniquely amenable topological group
such that $ \UMeas(\ru G) = Z $,
then $G$ is precompact.
\end{theorem}

\noindent
{\bf Proof.} If $ \UMeas(\ru G) = Z $ and $G$ is uniquely amenable then
there exists a left-invariant mean $ \mu \in \UMeas (\ru G) $
by Lemma~\ref{lemma:uniqmean}.
By Theorem~2 in~\cite{Pachl1982}, $G$ is precompact.
\qed

The next theorem shows that a result proved by Granirer~(\cite{Granirer1967}, Th.~4)
for a class of countably additive measures holds more generally for uniform measures.
A similar result for locally compact semigroups was proved by Wong
(\cite{Wong1972}, 3.1.(d)).

For $\mu\iin\UMeas(\ru G)$ and $f\iin\Ub(\ru G)$,
let $L_\mu f : G \rightarrow \real$ be the function
$\sect{y}\mu(\sect{x}f(xy))$.

\begin{theorem}
    \label{th:invmeans}
If $ \nu \in \Meas(\ru G) $ is a left-invariant mean,
$ \mu \in \UMeas( \ru G ) $ and $f\iin\Ub(\ru G)$
then
$ \nu(L_\mu f) = \mu \conv \nu (f)= \mu(1) \, . \, \nu (f)$.
\end{theorem}

\noindent
{\bf Proof.}
As is noted above, $ \UMeas(\ru G) \subseteq Z $ and therefore
\[
\nu(L_\mu f)
= \nu ( \sect{y} \mu ( \sect{x} f ( xy ) ) )
= \mu ( \sect{x} \nu ( \sect{y} f ( xy ) ) )
= \mu \conv \nu (f) \; .
\]
Since $\nu$ is left-invariant, $\mu \conv \nu (f)= \mu(1) \, . \, \nu (f)$
for every $\mu\in\Mol(\ru G)$.
Since $\Mol(\ru G)$ is \wstar dense in $\Meas(\ru G)$ and the mapping
$\mu\mapsto \mu\conv\nu$ is \wstar continuous from $\Meas(\ru G)$ to itself,
it follows that $\mu \conv \nu (f)= \mu(1) \, . \, \nu (f)$ for every $\mu\in\Meas(\ru G)$.
\qed


\section{Completion and compactification}
    \label{section:compactification}

For any (Hausdorff) uniform space $X$,
consider the mapping $ \delta : x \mapsto \delta_x $
from $X$ to $ \Meas (X) $,
where $ \delta_x $ is the Dirac measure at $x$.
The mapping $ \delta $ is a uniform embedding of $X$ into $ \Meas(X) $
with the UEB uniformity.
In the following discussion we
identify $X$ with its image $ \delta(X) \subseteq \Meas(X) $.

Define $ \compactn{X} $ to be the \wstar closure of $ X $ in $ \Meas(X) $,
and make $ \compactn{X} $ into a topological space by endowing it
with the \wstar topology.
Define $ \compln{X} = \compactn{X} \cap \UMeas(X) $,
and make $ \compln{X} $ into a uniform space by endowing it with
the UEB uniformity.

Then $ \compactn{X} $ is compact and every uniformly continuous mapping
from $X$ to a compact space extends uniquely to $ \compactn{X} $.
Thus $ \compactn{X} $ is a {\em uniform compactification\/}
(or Samuel compactification) of~$X$
(\cite{Isbell1964}, 2.32).

\begin{lemma}
    \label{lemma:compactification}
Let $X$ be a uniform space.
The following properties of $ \mu \in \Meas(X) $ are equivalent:
\begin{description}
\item[{\it (i)}]
$ \mu \in \compactn{X} $.
\item[{\it (ii)}]
$ \mu \neq 0 $ and
$ \mu (f g ) = \mu(f) . \, \mu(g) $ for all $ f,g \in \Ub(X) $.
\end{description}
\end{lemma}

\noindent
{\bf Proof.} This is a special case of (C.32) in~\cite{Hewitt-Ross1979}.
Cf. also 1.9 in \cite{BJM1989}.
\qed

Clearly $ \compln{X} \subseteq \UMeasplus(X) $, and
$ \compln{X} $ is the \wstar closure of $ X $ in $ \UMeas(X) $.
By Theorem~\ref{th:umeas-topology} (part~\ref{th:umeas-topology:part8}),
$ \compln{X} $ is also the UEB closure of $ X $ in $ \UMeas(X) $.
By part~\ref{th:umeas-topology:part3} of the same theorem,
the space $ \compln{X} $ is complete.
Thus $ \compln{X} $ is a completion of $X$.

By Lemma~\ref{lemma:compactification}, $ \mu \in \compln{X} $ if and only if
$ \mu \in \UMeas(X) $ and $ \mu $ is multiplicative.
This and related results are discussed by
Buchwalter and Pupier~\cite{Buchwalter-Pupier1971}.

From the definition of $ \compactn{X} $ and $ \compln{X} $
we have $ X \subseteq \compln{X} \subseteq \compactn{X} \subseteq \Meas(X) $.
The uniform compactification and the completion of $X$ may be constructed in other
ways, but they are unique up to a natural isomorphism.

Now consider the uniform compactification $\compactn{\ru G}$
and the completion $\compln{\ru G}$ of a topological group $G$
with its right uniformity.
Since $\compactn{\ru G} \subseteq \Meas ( \ru G ) $,
the convolution $ \mu \conv \nu $ is defined for
$ \mu, \nu \in \compactn{\ru G} $.

\begin{lemma}
Let $G$ be a topological group.
\begin{enumerate}
\item
If $ \mu, \nu \in \compactn{\ru G} $ then $ \mu \conv \nu \in \compactn{\ru G}$.
\item
If $ \mu, \nu \in \compln{\ru G} $ then $ \mu \conv \nu \in \compln{\ru G}$.
\end{enumerate}
\end{lemma}

\noindent
{\bf Proof.}
If $\mu$ and $\nu$ are multiplicative then so is $ \mu \conv \nu $.
Thus part~1 follows from Lemma~\ref{lemma:compactification}.

Part~2 follows from part~1 and Theorem~\ref{th:conv-closure}.
\qed

In the terminology of~\cite{BJM1989},
$ \compactn{ \ru G } $ with the $ \conv $ operation
is a semigroup compactification of $G$.
In fact,
$ \compactn{ \ru G } = G^{\cal L C}$ is the canonical $\cal L C$-compactification of $G$
(\cite{BJM1989}, 4.4).

As in~\cite{BJM1989}, define
\[
\Lambda ( \compactn{ \ru G } ) \; = \;
\{ \; x \in \compactn{ \ru G } \; | \;
\mbox{\rm the mapping} \;\; y \mapsto x \conv y  \;\;
\mbox{\rm is \wstar continuous on} \; \compactn{ \ru G } \; \} \; .
\]

From part~\ref{th:convolution-topology:part5} of
Theorem~\ref{th:convolution-topology}
we get $ \compln{ \ru G } \subseteq \Lambda ( \compactn{ \ru G } ) $
--- cf. Prop.~4.11 in~\cite{Ferri-Neufang2006}.
As with the question $ \UMeas(\ru G) \stackrel{?}{=} Z $ in the previous section,
one could ask whether
$ \compln{ \ru G } = \Lambda ( \compactn{ \ru G } ) $
for every topological group $G$.

If $G$ locally compact then $ \ru G $ is complete and therefore
$ \compln{ \ru G } = G $.
Lau and Pym~\cite{Lau-Pym1995} proved that $ G = \Lambda ( \compactn{ \ru G } ) $
for every locally compact group $ G $.
Ferri and Neufang~\cite{Ferri-Neufang2006} proved that
$ \compln{ \ru G } = \Lambda ( \compactn{ \ru G } ) $ for all
$\omega$-bounded groups $G$.


\section{Generalized convolution}
    \label{section:conclude}

In section~\ref{section:convolution},
the inverse group operation and the topology on $G$ are used only to establish
the conclusion of Lemma~\ref{lemma:rightu}.
That suggests the following definition.
$S$ is a {\em (right) semiuniform semigroup\/} if $S$ is a semigroup and
a (Hausdorff) uniform space,
and the semigroup operation in $S$ is uniformly continuous
from the semiuniform product $ S \semiu S $ to $S$.

By Lemma~\ref{lemma:rightu},
any topological group with its right uniformity is a semiuniform semigroup.
Uniform semigroups as defined by Marxen~\cite{Marxen1973} are semiuniform semigroups.

For any semiuniform semigroup $S$ we may define the convolution operation on
$ \Meas(S) $ as in section~\ref{section:convolution}.
Namely,
$
\mu \conv \nu = \Meas ( m_S ) ( \mu \measprod \nu )
$
where $ m_S $ is the binary semigroup operation in $S$.
All the results in section~\ref{section:convolution},
with the exception of Theorem~\ref{th:metrizableexample},
hold for semiuniform semigroups
in place of topological groups, with the same proofs.

Now consider the following generalization of the approach outlined on the last page
of Csisz\'{a}r's paper~\cite{Csiszar1971}.
Let $S$ be a semiuniform semigroup acting on a uniform space $Y$, in the sense that
there is a uniformly continuous mapping $ \alpha: S \semiu Y \rightarrow Y $ such that
$ \alpha ( s , \alpha ( s' , y ) ) = \alpha ( s s' , y ) $
for all $ s, s' \in S $, $ y \in Y$.
Define the convolution operation from $ \Meas(S) \times \Meas(Y) $ to $ \Meas(Y) $
by $ \mu \conv \nu = \Meas( \alpha ) ( \mu \measprod \nu ) $.

Lemma~\ref{lemma:dproduct3} shows that the semigroup $ \Meas(S) $
acts on $ \Meas(Y) $ using the operation $\conv$,
and the results in section~\ref{section:directp} immediately yield continuity
properties of the $\conv$ operation analogous to those in
Theorem~\ref{th:convolution-topology}.


\end{document}